\documentclass[a4paper,12pt,reqno]{amsart}

\usepackage[mathscr]{eucal}
\usepackage{amscd}
\usepackage{amsfonts}
\usepackage{amsmath, amsthm, amssymb}
\usepackage{latexsym}
\usepackage[dvips]{graphics}

\usepackage{lipsum} 
\usepackage{cite}
\usepackage[sc]{mathpazo} 
\usepackage[T1]{fontenc} 
\usepackage{microtype} 

\usepackage[hmarginratio=1:1,top=30mm]{geometry} 
\usepackage{hyperref} 
\usepackage{abstract} 

\usepackage{titlesec} 
\renewcommand\thesection{\Roman{section}} 
\titleformat{\section}[block]{\large\scshape\centering}{\thesection.}{1em}{} 
\renewcommand{\thesection}{\arabic{section}}

\newtheorem{thm}{Theorem}[section]

\numberwithin{equation}{section}

\newtheorem{prop}{Proposition}[section]
\newtheorem{cor}{Corollary}[section]

\theoremstyle{plain}

\begin{document}

\thispagestyle{empty}



\title{Multiplicative generalized Jordan $n$-derivations of unital rings with idempotents \\     
}


\author{Mohammad Ashraf, Mohammad Afajal Ansari and Md Shamim Akhter}

\address{Mohammad Ashraf, Department of Mathematics,
Aligarh Muslim University,
Aligarh-202002 India}
\email{\href{mailto:mashraf80@hotmail.com}{mashraf80@hotmail.com}}

\address{Mohammad Afajal Ansari, Department of Mathematics, Aligarh Muslim University,
Aligarh-202002 India}
\email{\href{mailto:afzalgh1786@gmail.com}{afzalgh1786@gmail.com}}

\address{Md Shamim Akhter, Department of Mathematics,
Aligarh Muslim University,
Aligarh-202002 India}
\email{\href{mailto:akhter2805@gmail.com}{akhter2805@gmail.com}}

\maketitle



\begin{abstract}
Let $\mathfrak{A}$ be a unital ring with a nontrivial idempotent. In this paper, it is shown that under certain conditions every multiplicative generalized Jordan $n$-derivation $\Delta:\mathfrak{A}\rightarrow\mathfrak{A}$ is additive. More precisely, it is proved that $\Delta$ is of the form $\Delta(t)=\mu t+\delta(t),$ where $\mu\in\mathcal{Z}(\mathfrak{A})$ and $\delta:\mathfrak{A}\rightarrow\mathfrak{A}$ is a Jordan $n$-derivation. The main result is then applied to some classical examples of unital rings with nontrivial idempotents such as triangular rings, matrix rings, prime rings, nest algebras, standard operator algebras and von Neumann algebras.
\end{abstract}

\vspace{.40cm}



\vspace{.40cm}

\section{Introduction}
Let $\mathfrak{A}$ be an associative ring with identity. An additive mapping $\delta:\mathfrak{A}\rightarrow \mathfrak{A}$ is called a \textit{derivation} if $\delta(t_1t_2)=\delta(t_1)t_2+t_1\delta(t_2)$ holds for all $t_1, t_2\in \mathfrak{A}$. An additive mapping $\delta:\mathfrak{A}\rightarrow \mathfrak{A}$ is said to be a \textit{Jordan derivation} (resp. \textit{Jordan triple derivation}) if $\delta(t_1\circ t_2) = \delta(t_1)\circ t_2 + t_1\circ \delta(t_2)$ (resp.  $\delta((t_1\circ t_2) \circ t_3)=(\delta(t_1)\circ t_2)\circ t_3+(t_1\circ \delta(t_2))\circ t_3+ (t_1\circ t_2)\circ \delta(t_3)$) holds for all $t_1, t_2, t_3\in \mathfrak{A},$ where $t_1\circ t_2=t_1t_2+t_2t_1$ is the usual Jordan product. Note that every derivation is a Jordan derivation and every Jordan derivation is a Jordan triple derivation but the converse statements are not true in general. The natural question that arises in such context is that under what additional conditions on a given algebra every Jordan (triple) derivation is a derivation. Such type of problems have been studied extensively by many authors (see \cite{ap16,aja19,h57,c75,b89,zy06,Gh07,ll07,fi08,bs12} and references therein).

Further, the notions of Jordan derivation, Jordan triple derivation can be extended to a more general class of mappings. Set $q_1(t_1)=  t_1$ and
$$q_n(t_1, t_2,\cdots,t_n)=q_{n-1}(t_1,t_2,\cdots,t_{n-1})\circ t_n \ \ \text{for all integer} \ \ n\geq 2.$$
An additive mapping $\delta:\mathfrak{A}\rightarrow\mathfrak{A}$ is called a \textit{Jordan ${n}$-derivation} if
\begin{equation}\label{eqn1.1}
  \delta(q_{n}(t_{1},t_{2},\cdots,t_{n}))=\sum\limits_{{i=1}}^{{n}}q_{n}(t_{1},t_
{2},\cdots,t_{i-1},\delta(t_{i}),t_{i+1},\cdots,t_{n})
\end{equation}
for all $t_{1},t_{2},\cdots,t_{n}\in\mathfrak{A}.$ One can easily see that a Jordan ${2}$-derivation is the usual Jordan derivation and a Jordan ${3}$-derivation is a Jordan triple derivation. An additive mapping $\Delta:\mathfrak{A}\rightarrow\mathfrak{A}$ is said to be a \textit{generalized Jordan  ${n}$-derivation} if there exists a Jordan ${n}$-derivation $\delta:\mathfrak{A}\rightarrow\mathfrak{A}$ such that
\begin{eqnarray}\label{eqn1.2}
\nonumber  \Delta(q_{n}(t_{1},t_{2},\cdots,t_{n})) &=& q_n(\Delta(t_1),t_2,\cdots,t_{n-1},t_n) \\
   &&+ \sum\limits_{{i=2}}^{{n}}q_{n}(t_{1},t_{2},\cdots,t_{i-1},\delta(t_{i}),t_{i+1},\cdots,t_{n})
\end{eqnarray}
for all $t_{1},t_{2},\cdots,t_{n}\in\mathfrak{A}.$ In case the mappings $\delta$ and $\Delta$ are not necessarily additive in the above definitions, $\Delta$ is called a multiplicative generalized Jordan ${n}$-derivation with associated multiplicative Jordan ${n}$-derivation $\delta.$ Note that any multiplicative Jordan $n$-derivation and the mapping $t\mapsto \mu t,$ where $\mu\in\mathcal{Z}(\mathfrak{A}),$ are examples of multiplicative generalized Jordan $n$-derivation.

A lot of work has been done on the additivity of mappings on various rings and algebras over the past decade. In the year 1969, Martindale III \cite{m69} proved  that under certain assumptions any multiplicative  isomorphism from a prime ring containing a nontrivial idempotent onto an arbitrary ring is additive. Daif \cite{d91} proved that every multiplicative derivation of a $2$-torsion free prime ring containing a nontrivial idempotent is additive. Inspired by these results many researchers obtained similar results in various rings and algebras (see \cite{l10,jl12,lf13} and references therein).  Recently, Qi et al. \cite{qgz19} introduced the notion of a multiplicative Jordan $n$-derivation and discussed additivity of multiplicative Jordan $n$-derivations of unital rings containing nontrivial idempotents.

In this article we study the additivity and structure of multiplicative generalized Jordan $n$-derivations of unital rings with idempotents. In fact, we prove that under certain conditions every multiplicative generalized Jordan $n$-derivation $\Delta$ of a unital ring $\mathfrak{A}$ is additive. More precisely, it is shown that $\Delta$ is of the form $\Delta(t)=\mu t+\delta(t)$ for all $t\in\mathfrak{A}$, where $\mu\in\mathcal{Z}(\mathfrak{A})$ and $\delta$ is a Jordan $n$-derivation of $\mathfrak{A}$ (Theorem \ref{thm3.1}). As a consequence of the main result, multiplicative generalized Jordan $n$-derivations of triangular rings, nest algebras, matrix rings, prime rings, standard operator algebras and von Neumann algebras will be determined. The main motivation of this article comes from the papers \cite{qgz19,B18,lb11}. Li and Benkovi\v c in \cite{lb11} studied Jordan generalized derivation (i.e. linear generalized Jordan $2$-derivation) of triangular algebras and proved that any Jordan generalized derivation of a triangular algebra is a generalized derivation. Recently, Benkovi\v{c} \cite{B18} described generalized Lie $n$-derivations of unital algebras with nontrivial idempotents and showed that under certain mild restrictions every generalized Lie $n$-derivation $\Delta$ on a unital algebra $\mathfrak{A}$ can be written as $\Delta(t)=\mu t+\delta(t)$ for all $t\in\mathfrak{A}$, where $\mu \in\mathcal{Z}(\mathfrak{A})$ and $\delta$ is a Lie $n$-derivation of $\mathfrak{A}.$

Let $\Delta:\mathfrak{A}\rightarrow\mathfrak{A}$ be a multiplicative generalized Jordan $n$-derivation with associated multiplicative Jordan $n$-derivation $\delta:\mathfrak{A}\rightarrow\mathfrak{A}.$ Set $\Phi=\Delta-\delta.$ Then it follows from $(\ref{eqn1.1})$ and $(\ref{eqn1.2})$ that $\Phi$ satisfies
\begin{eqnarray}\label{eqn1.3}
  \Phi(q_{n}(t_{1},t_{2},\cdots,t_{n})) &=& q_{n}(\Phi(t_{1}),t_2, \cdots,t_{n})
\end{eqnarray}
for all $t_{1},t_{2},\cdots,t_{n}\in\mathfrak{A}.$ Therefore, to characterize multiplicative generalized Jordan $n$-derivations, it suffices to consider a mapping satisfying $(\ref{eqn1.3}).$ Any mapping $\Phi:\mathfrak{A}\rightarrow\mathfrak{A}$ satisfying $(\ref{eqn1.3})$ is called a \textit{multiplicative Jordan $n$-centralizer}. In particular, multiplicative Jordan $2$-centralizer is a mapping $\Phi:\mathfrak{A}\rightarrow\mathfrak{A}$ that satisfies $\Phi(t_1\circ t_2)=\Phi(t_1)\circ t_2$ for all $t_1, t_2\in\mathfrak{A}.$ Recall that an additive mapping $\Phi:\mathfrak{A}\rightarrow \mathfrak{A}$ is called a \textit{centralizer} if $\Phi(t_1t_2)=\Phi(t_1)t_2=t_1\Phi(t_2)$ for all $t_1, t_2\in \mathfrak{A}.$  If $\mathfrak{A}$ is a unital ring, then $\Phi$ is a centralizer if and only if $\Phi(t)=\mu t$ for some $\mu\in\mathcal{Z}(\mathfrak{A}).$ In 1991, Zalar \cite{Z91} introduced the notion of Jordan centralizers and proved that every Jordan centralizer on a $2$-torsion free semiprime ring is a centralizer. Vukman and Kosi-Ulbl \cite{j99,j01,jk06} extensively studied centralizers mainly on prime rings and semiprime rings. Motivated by these results, we show that under certain restrictions every multiplicative Jordan $n$-centralizer $\Phi$ on a unital ring $\mathfrak{A}$ is of the form $\Phi(t)=\mu t $ for all $t\in\mathfrak{A}$, where $\mu\in\mathcal{Z}(\mathfrak{A})$ (Proposition \ref{prop3.1}). Applying this result, we obtain the main result of the paper (Theorem \ref{thm3.1}).

\section{Preliminaries}
Let $\mathfrak{A}$ be a unital ring with a nontrivial idempotent $e_1,$ and write $e_2=1-e_1.$ Then $\mathfrak{A}$ can be represented as $\mathfrak{A}=e_1\mathfrak{A}e_1+e_1\mathfrak{A}e_2+e_2\mathfrak{A}e_1+e_2\mathfrak{A}e_2,$ where $e_1\mathfrak{A}e_1$ and $e_2\mathfrak{A}e_2$ are subrings of $\mathfrak{A}$ with identity elements $e_1$ and $e_2,$ respectively, $e_1\mathfrak{A}e_2$ is an $(e_1\mathfrak{A}e_1,e_2\mathfrak{A}e_2)$-bimodule and $e_2\mathfrak{A}e_1$ is an $(e_2\mathfrak{A}e_2,e_1\mathfrak{A}e_1)$-bimodule. Throughout the paper, we assume that $\mathfrak{A}$ is a $2$-torsion free unital ring with a nontrivial idempotent $e$ satisfying the following conditions:
\[\begin{cases}
 e_1te_1\cdot e_1\mathfrak{A}e_2=\{0\}=e_2\mathfrak{A}e_1\cdot e_1te_1~~~~ \mbox{implies}~~~~ e_1te_1=0\\
e_1\mathfrak{A}e_2\cdot e_2te_2=\{0\}=e_2te_2\cdot e_2\mathfrak{A}e_1~~~~ \mbox{implies}~~~~ e_2te_2=0 \tag{$\spadesuit$}
 \end{cases} \]
for all $t\in\mathfrak{A}.$ Some standard examples of unital rings with nontrivial idempotents satisfying $(\spadesuit)$ are triangular rings, matrix rings, unital prime rings with nontrivial idempotents and nest algebras. To simplify the calculations, we will use the following notations: $t_{11}=e_1te_1\in e_1\mathfrak{A}e_1=\mathfrak{A}_{11}$, $t_{12}=e_1te_2\in e_1\mathfrak{A}e_2=\mathfrak{A}_{12},$ $t_{21}=e_2te_1\in e_2\mathfrak{A}e_1=\mathfrak{A}_{21}$ and $t_{22}=e_2te_2\in e_2\mathfrak{A}e_2=\mathfrak{A}_{22}.$ Thus each element $t \in \mathfrak{A}$ can be written as $t=t_{11}+t_{12}+t_{21}+t_{22}.$

In view of \cite[Proposition 2.1]{B15}, the center of $\mathfrak{A}$ is given by
$$\mathcal{Z}(\mathfrak{A})=\bigg\{t_{11}+t_{22}~~\vline~~
\begin{array}{c}
t_{11}\in e_1\mathfrak{A}e_1,~ t_{22}\in e_2\mathfrak{A}e_2,~  t_{11}t_{12} = t_{12}t_{22},   \\
t_{21}t_{11} = t_{22}t_{21}~\text{for all}~~t_{12}\in \mathfrak{A}_{12},~ t_{21}\in \mathfrak{A}_{21}
                                                          \end{array}
\bigg\}.$$
Furthermore, there exists a unique algebra isomorphism $\xi:\mathcal{Z}(\mathfrak{A})e_1\rightarrow \mathcal{Z}(\mathfrak{A})e_2$ such that $t_{11}t_{12} = t_{12}\xi(t_{11})$ and $t_{21}t_{11} = \xi(t_{11})t_{21}$ for all $t_{12}\in e_1\mathfrak{A}e_2$, $t_{21}\in e_2\mathfrak{A}e_1$ and for any $t_{11}\in\mathcal{Z}(\mathfrak{A})e_1.$

Suppose that $\mathcal{R}$ is a ring such that for each $t \in \mathcal{R}$
\begin{eqnarray}\label{eqn2.1}
[t,\mathcal{R}] \subseteq \mathcal{Z}(\mathcal{R})\Longrightarrow t \in \mathcal{Z}(\mathcal{R}),
\end{eqnarray}
It is straightforward to prove that $\mathfrak{A}$ satisfies the condition (\ref{eqn2.1}). For some more examples of rings satisfying the condition (\ref{eqn2.1}), we refer the reader to \cite[Section 5]{be12}.

\section{Multiplicative generalized Jordan $n$-derivation}
We begin this section with the following proposition which characterizes multiplicative Jordan $n$-centralizers of unital rings with nontrivial idempotents satisfying $(\spadesuit).$

\begin{prop}\label{prop3.1}
Let $\mathfrak{A}$ be a $2$-torsion free unital ring with a nontrivial idempotent $e_1$ satisfying $(\spadesuit).$
Then every multiplicative Jordan $n$-centralizer $\Phi:\mathfrak{A}\rightarrow\mathfrak{A}$ is of the form
$\Phi(t)= \mu t$ for all $t\in \mathfrak{A},$ where $\mu\in\mathcal{Z}(\mathfrak{A}).$
\end{prop}

\begin{proof}
Let $\Phi:\mathfrak{A}\rightarrow\mathfrak{A}$ be a multiplicative Jordan $n$-centralizer. Then
\begin{eqnarray}\label{eqn4.3.2}
  \Phi(q_{n}(t_{1},t_{2},\cdots,t_{n}))=q_{n}(\Phi(t_{1}),t_2, \cdots,t_{n})
\end{eqnarray}
for all $t_{1},t_{2},\cdots,t_{n}\in\mathfrak{A}.$ Since Jordan product is commutative, \eqref{eqn4.3.2} yields
\begin{eqnarray}\label{eqn4.3.3}
\nonumber  \Phi(q_{n}(t_{1},t_{2},\cdots,t_{n})) &=& \Phi(q_{n}(t_{2},t_{1},\cdots,t_{n}))\\
\nonumber  &=&q_{n}(\Phi(t_{2}),t_1, \cdots,t_{n})\\
  &=&q_{n}(t_1,\Phi(t_{2}), \cdots,t_{n})
\end{eqnarray}
for all $t_{1},t_{2},\cdots,t_{n}\in\mathfrak{A}.$
Combining \eqref{eqn4.3.2} and \eqref{eqn4.3.3}, we get
\begin{eqnarray}\label{eqn4.3.4}
  q_{n}(\Phi(t_{1}),t_{2},\cdots,t_{n}))=q_{n}(t_1,\Phi(t_{2}), \cdots,t_{n})
\end{eqnarray}
for all $t_{1},t_{2},\cdots,t_{n}\in\mathfrak{A}.$
Let us denote $\Phi(1)=\mu\in\mathfrak{A}$. If we insert $t_1 = t$ and $t_2 =\ldots= t_n = 1$ into \eqref{eqn4.3.4},
note that $\mathfrak{A}$ is $2$-torsion free, we get
\begin{equation}\label{eqn4.3.5}
 2\Phi(t) = \mu \circ t
\end{equation}
for all $t\in\mathfrak{A}$. Next, if we insert $t_1 = t,$ $t_2 = t^{\prime}$ and $t_3=\ldots=t_n = 1$ into \eqref{eqn4.3.2},
we have
\begin{equation*}
  \Phi(2^{n-2}(t\circ t^{\prime}))=2^{n-2}(\Phi(t)\circ t^{\prime})
\end{equation*}
for all $t,t^{\prime}\in\mathfrak{A},$ which can be further rewritten as
$$2\Phi(2^{n-2}(t\circ t^{\prime}))=2^{n-2}(2\Phi(t)\circ t^{\prime}).$$
In view of \eqref{eqn4.3.5}, the above equation gives
\begin{eqnarray*}
  \mu \circ(2^{n-2}(t\circ t^{\prime})) &=& 2^{n-2}((\mu \circ t)\circ t^{\prime})\\
  \mu \circ (t\circ t^{\prime})&=&(\mu \circ t)\circ t^{\prime}
\end{eqnarray*}
for all $t,t^{\prime}\in\mathfrak{A}.$
Therefore, $[[\mu, t],t^{\prime}]=\mu \circ (t\circ t^{\prime})-(\mu \circ t)\circ t^{\prime}=0$ for all $t,t^{\prime}\in\mathfrak{A}.$
Since $\mathfrak{A}$ satisfies the condition \eqref{eqn2.1},
we conclude that $\mu\in\mathcal{Z}(\mathfrak{A}).$ Hence, it follows from \eqref{eqn4.3.5} that
$\Phi(t)= \mu t$ for all $t\in \mathfrak{A},$ where $\mu\in\mathcal{Z}(\mathfrak{A}).$
\end{proof}

Applying Proposition \ref{prop3.1}, we now obtain the main result of the paper which characterizes multiplicative generalized Jordan  $n$-derivations of unital rings with nontrivial idempotents satisfying $(\spadesuit).$
\begin{thm}\label{thm3.1}
Let $\mathfrak{A}$ be a $2$-torsion free unital ring with a nontrivial idempotent $e_1$ satisfying $(\spadesuit).$ Then every multiplicative generalized Jordan $n$-derivation $\Delta:\mathfrak{A}\rightarrow\mathfrak{A}$ is additive. More precisely, $\Delta$ is of the form $\Delta(t)=\mu t+\delta(t)$ for all $t\in \mathfrak{A},$ where $\mu\in\mathcal{Z}(\mathfrak{A})$ and $\delta:\mathfrak{A}\rightarrow \mathfrak{A}$ is a Jordan $n$-derivation.
\end{thm}

\begin{proof}
Let us assume that $\Delta:\mathfrak{A}\rightarrow\mathfrak{A}$ be a multiplicative generalized Jordan $n$-derivation with associated multiplicative Jordan $n$-derivation $\delta:\mathfrak{A}\rightarrow\mathfrak{A}.$ In view of \cite[Theorem 2.1]{qgz19}, $\delta$ is a Jordan $n$-derivation. Set $\Phi=\Delta-\delta.$ By the definitions of $\Delta,\delta,$ it follows that $\Phi:\mathfrak{A}\rightarrow\mathfrak{A}$ is a multiplicative Jordan $n$-centralizer. Hence it follows from Proposition \ref{prop3.1} that there exists a $\mu\in\mathcal{Z}(\mathfrak{A})$ such that $\Phi(t)= \mu t$ for all $t\in \mathfrak{A}.$ It is easy to see that $\Phi$ is additive and hence $\Delta$ is additive. Moreover, $\Delta$ is of the form $\Delta(t)=\mu t+\delta(t)$ for all $t\in \mathfrak{A},$  as desired.
\end{proof}

Recall that a Jordan derivation $\phi:\mathfrak{A}\rightarrow \mathfrak{A}$ is said to be a singular Jordan derivation if
$$
\phi(e_1\mathfrak{A}e_1)=\{0\}, \ \  \phi(e_1\mathfrak{A}e_2)\subseteq e_2\mathfrak{A}e_1, \ \  \phi(e_2\mathfrak{A}e_1)\subseteq e_1\mathfrak{A}e_2, \ \ \phi(e_2\mathfrak{A}e_2)=\{0\}.
$$
Let $\mathfrak{A}$ be a $2$-torsion free unital ring with a nontrivial idempotent $e_1$ satisfying $(\spadesuit).$
It is shown in \cite[Theorem 3.3]{qgz19} that every Jordan $n$-derivation $\delta:\mathfrak{A}\rightarrow\mathfrak{A}$
is of the form $\delta(u)=d(u)+\phi(u)$ for all $u\in \mathfrak{A},$ where $d:\mathfrak{A}\rightarrow \mathfrak{A}$
is a derivation and $\phi:\mathfrak{A}\rightarrow \mathfrak{A}$ is a singular Jordan derivation.
Combining the above result with Theorem \ref{thm3.1}, we immediately obtain the following result:

\begin{cor}
Let $\mathfrak{A}$ be a $2$-torsion free unital ring with a nontrivial idempotent $e_1$ satisfying $(\spadesuit).$ Then every multiplicative generalized Jordan $n$-derivation $\Delta:\mathfrak{A}\rightarrow\mathfrak{A}$ is of the form $\Delta(t)=\mu t+d(t)+\phi(t)$ for all $t\in \mathfrak{A},$ where $\mu\in\mathcal{Z}(\mathfrak{A}),$ $d:\mathfrak{A}\rightarrow \mathfrak{A}$ is a derivation and $\phi:\mathfrak{A}\rightarrow \mathfrak{A}$ is a singular Jordan derivation.
\end{cor}

Let $\mathfrak{A}$ be a unital ring with a nontrivial idempotent $e_1$ such that
$e_1\mathfrak{A}e_2\cdot e_2\mathfrak{A}e_1=\{0\}=e_2\mathfrak{A}e_1\cdot e_1\mathfrak{A}e_2.$
It is shown in \cite[Remark 3.2]{bs12} that every singular Jordan derivation of $\mathfrak{A}$ is an antiderivarion;
that is, an additive mapping $\phi:\mathfrak{A}\rightarrow \mathfrak{A}$ satisfying $\phi(xy)=\phi(y)x+y\phi(x)$
for all $x,y\in\mathfrak{A}.$
If the bimodule $e_1\mathfrak{A}e_2$ is faithful as a
left $e_1\mathfrak{A}e_1$-module and also as a right $e_2\mathfrak{A}e_2$-module, then $\mathfrak{A}$ satisfies $(\spadesuit)$
and we obtain:
\begin{cor}
Let $\mathfrak{A}$ be a unital ring with a nontrivial idempotent $e_1$ such that
$e_1\mathfrak{A}e_2\cdot e_2\mathfrak{A}e_1=\{0\}=e_2\mathfrak{A}e_1\cdot e_1\mathfrak{A}e_2$ and the bimodule $e_1\mathfrak{A}e_2$ is faithful as a left $e_1\mathfrak{A}e_1$-module and also as a right $e_2\mathfrak{A}e_2$-module.
Then every multiplicative generalized Jordan $n$-derivation $\Delta:\mathfrak{A}\rightarrow\mathfrak{A}$ is of the form $\Delta(t)=\mu t+d(t)+\phi(t)$ for all $t\in \mathfrak{A},$ where $\mu\in\mathcal{Z}(\mathfrak{A}),$ $d:\mathfrak{A}\rightarrow \mathfrak{A}$ is a derivation and $\phi:\mathfrak{A}\rightarrow \mathfrak{A}$ is an antiderivation.
\end{cor}

\section{Applications}
In this section, we apply Theorem \ref{thm3.1} to certain classes of unital rings satisfying $(\spadesuit)$ such as triangular rings, prime rings containing nontrivial idempotents, matrix rings and nest algebras.

\subsection*{Triangular Rings:} Let $\mathfrak{A}$ be a unital ring with a nontrivial idempotent $e_1$ such that $e_1\mathfrak{A}e_2$ is a faithful $(e_1\mathfrak{A}e_1,e_2\mathfrak{A}e_2)$-bimodule and $e_2\mathfrak{A}e_1=\{0\}.$ Then $\mathfrak{A}=e_1\mathfrak{A}e_1+e_1\mathfrak{A}e_2+e_2\mathfrak{A}e_2$ is a triangular ring. One can easily see that $\mathfrak{A}$ satisfies the assumption $(\spadesuit).$ Therefore,  in view of  Theorem \ref{thm3.1} and \cite[Corollary 3.5]{qgz19}, we obtain the following result:

\begin{cor}\label{cor4.1}
Let $\mathfrak{A}=e_1\mathfrak{A}e_1+e_1\mathfrak{A}e_2+e_2\mathfrak{A}e_2$ be a $2$-torsion free triangular ring.
Then every multiplicative generalized  Jordan $n$-derivation $\Delta:\mathfrak{A}\rightarrow\mathfrak{A}$ can be expressed as $\Delta(t)=\mu t+\delta(t)$ for all $t\in \mathfrak{A},$ where $\mu\in\mathcal{Z}(\mathfrak{A})$ and $\delta:\mathfrak{A}\rightarrow \mathfrak{A}$ is a derivation.
\end{cor}

Some standard examples of triangular rings are upper triangular matrix rings, block upper triangular matrix rings and nest algebras
(see \cite{be12} for details). Hence, applying Corollary \ref{cor4.1}, we obtain the following results:

\begin{cor}{\label{cor4.2}}
Let $\mathfrak{A}$ be a $2$-torsion free unital ring and $T_r({\mathfrak{A}})(r\geq 2)$ be a upper triangular matrix ring. Then every multiplicative generalized  Jordan $n$-derivation $\Phi:T_r({\mathfrak{A}})\rightarrow T_r({\mathfrak{A}})$ can be expressed as $\Delta(t)=\mu t+\delta(t)$ for all $t\in T_r({\mathfrak{A}}),$ where $\mu\in\mathcal{Z}(\mathfrak{A})$ and $\delta:T_r({\mathfrak{A}})\rightarrow T_r({\mathfrak{A}})$ is a derivation.
\end{cor}

\begin{cor}{\label{cor4.3}}
Let $\mathfrak{A}$ be a $2$-torsion free unital ring and $B^{\overline{k}}_r({\mathfrak{A}})(r\geq 2)$ be a block upper triangular matrix ring with $B^{\overline{k}}_r({\mathfrak{A}})\neq M_r(\mathfrak{A}).$ Then every multiplicative generalized Jordan $n$-derivation $\Delta:B^{\overline{k}}_r({\mathfrak{A}})\rightarrow B^{\overline{k}}_r({\mathfrak{A}})$ can be expressed as $\Delta(t)=\mu t+\delta(t)$ for all $t\in B^{\overline{k}}_r({\mathfrak{A}}),$ where $\mu\in\mathcal{Z}(\mathfrak{A})$ and $\delta:B^{\overline{k}}_r({\mathfrak{A}})\rightarrow B^{\overline{k}}_r({\mathfrak{A}})$ is a derivation.
\end{cor}

\begin{cor}{\label{cor4.4}}
Let $\mathcal{N}$ be a nontrivial nest of a complex Hilbert space $\bf {H}$ and $\mathcal{T}(\mathcal{N})$ be a nest algebra. Then every multiplicative generalized Jordan $n$-derivation $\Delta:\mathcal{T}(\mathcal{N})\rightarrow \mathcal{T}(\mathcal{N})$ can be expressed as $\Delta(t)=\mu t+\delta(t)$ for all $t\in \mathcal{T}(\mathcal{N}),$ where $\mu\in\mathbb{C}1$ and $\delta:\mathcal{T}(\mathcal{N})\rightarrow \mathcal{T}(\mathcal{N})$ is a derivation.
\end{cor}

\subsection*{Matrix Rings:} Let $M_{r}(\mathfrak{A})(r\geq 2)$ be a matrix ring over a unital ring $\mathfrak{A}.$ For each $1\leq k\leq r-1,$ the matrix ring $M_{r}(\mathfrak{A})$ can be written as
\begin{eqnarray*}
M_{r}(\mathfrak{A}) &=& \left(
                       \begin{array}{cc}
                         \mathfrak{A} & M_{1\times {(r-1)}}(\mathfrak{A}) \\
                          M_{{(r-1)}\times{1}}(\mathfrak{A}) &  M_{(r-1)}(\mathfrak{A}) \\
                       \end{array}
                     \right).
\end{eqnarray*}
Set $e_1=e_{11},$ the matrix unit and $e_2=1-e_1.$ Then it is easy to see that $M_{r}(\mathfrak{A})$ satisfies $(\spadesuit).$ Moreover, $e_1M_{r}(\mathfrak{A})e_1$ is isomorphic to $\mathfrak{A},$ $e_2M_{r}(\mathfrak{A})e_2$ is isomorphic to $M_{(r-1)}(\mathfrak{A}),$ $e_1M_{r}(\mathfrak{A})e_2$ is faithful $(e_1M_{r}(\mathfrak{A})e_1, e_2M_{r}(\mathfrak{A})e_2)$-bimodule isomorphic to $M_{{(r-1)}\times{1}}(\mathfrak{A}).$ Therefore, by Theorem \ref{thm3.1}, \cite[Theorem 3.3]{qgz19} and \cite[Corollary 2]{Gh07}, the following result holds:

\begin{cor}{\label{cor4.5}}
Let $\mathfrak{A}$ be a $2$-torsion free unital ring and $M_r(\mathfrak{A})$ be a matrix ring with $r\geq 2.$ Then every multiplicative  generalized Jordan $n$-derivation $\Delta:M_r(\mathfrak{A})\rightarrow M_r(\mathfrak{A})$ can be expressed as $\Delta(t)=\mu t+\delta(t)$ for all $t\in M_r(\mathfrak{A}),$ where $\mu\in\mathcal{Z}(\mathfrak{A})$ and $\delta:M_r(\mathfrak{A})\rightarrow M_r(\mathfrak{A})$ is a derivation.
\end{cor}

\subsection*{Prime Rings:} Recall that a ring $\mathfrak{A}$ is a called a prime ring if, for any $u, v\in\mathfrak{A},$ $u\mathfrak{A}v=\{0\}$ implies that $u=0$ or $v=0.$ It is easy to see that if $\mathfrak{A}$ has a nontrivial idempotent, then it satisfies the assumption $(\spadesuit).$
Therefore, as a consequence of Theorem \ref{thm3.1} and \cite[Corollary 3.10]{qgz19}, we obtain the following result:

\begin{cor}\label{cor4.6}
Let $\mathfrak{A}$ be a unital $2$-torsion free prime ring with a nontrivial idempotent. Then every multiplicative generalized Jordan  $n$-derivation $\Delta:\mathfrak{A}\rightarrow\mathfrak{A}$ can be expressed as $\Delta(t)=\mu t+\delta(t)$ for all $t\in \mathfrak{A},$ where $\mu\in\mathcal{Z}(\mathfrak{A})$ and $\delta:\mathfrak{A}\rightarrow \mathfrak{A}$ is a derivation.
\end{cor}

\subsection*{Standard Operator Algebras:} Let $\mathcal{B}(\mathcal{X})$ be the algebra of all bounded linear operators of a complex Banach space $\mathcal{X}$ of dimension greater than $1.$ A standard operator algebra $\mathcal{A}$ is a subalgebra of $\mathcal{B}(\mathcal{X})$ containing all bounded finite rank operators and identity operator. We know that $\mathcal{A}$ is a prime ring containing many nontrivial idempotents.
Therefore, Corollary \ref{cor4.6} yields the following:

\begin{cor}\label{cor4.7}
Let $\mathcal{A}$ be a standard operator algebra. Then every multiplicative generalized Jordan $n$-derivation $\Delta:\mathcal{A}\rightarrow\mathcal{A}$ can be expressed as $\Delta(t)=\mu t+\delta(t)$ for all $t\in \mathcal{A},$ where $\mu\in\mathcal{Z}(\mathcal{A})$ and $\delta:\mathcal{A}\rightarrow \mathcal{A}$ is a derivation.
\end{cor}

\subsection*{von Neumann Algebras:} Let $\mathcal{B}(\mathcal{H})$ be the algebra of all bounded linear operators of a complex Hilbert space $\mathcal{H}.$ A von Neumann algebra $\mathcal{M}$ is a $C^{\ast}$-subalgebra of $\mathcal{B}(\mathcal{H})$ which is closed in strong operator topology and contains the identity operator. If $\mathcal{M}$ has no central summands of type $I_1,$ then $\mathcal{M}$ satisfies $(\spadesuit)$ (see \cite[Section 3.5]{qgz19} for details).
Therefore, as a direct consequence of Theorem \ref{thm3.1} and \cite[Corollary 3.12]{qgz19}, we have the following result:
\begin{cor}\label{cor4.8}
Let $\mathcal{M}$ be a von Neumann algebra with no central summands of type $I_1.$ Then every multiplicative generalized Jordan $n$-derivation $\Delta:\mathcal{M}\rightarrow\mathcal{M}$ can be expressed as $\Delta(t)=\mu t+\delta(t)$ for all $t\in \mathcal{M},$ where $\mu\in\mathcal{Z}(\mathcal{M})$ and $\delta:\mathcal{M}\rightarrow \mathcal{M}$ is a derivation.
\end{cor}

\noindent
{\bf Acknowledgements}: The second author is partially supported by a research grant from DST (No. DST/INSPIRE/03/2017/IF170834).\\
\noindent
{\bf Conflict of Interests}: The authors have no conflicts of interest to declare that are relevant to the content of this article.\\
\noindent
{\bf Data Availability Statement}: This manuscript has no associate data.

\begin{center}

\end{center}

\end{document}